\documentclass[a4paper,12pt]{article}
\usepackage{comment}
\usepackage{cite}
\usepackage{amsmath}
\usepackage{amssymb}
\usepackage{amsfonts}
\usepackage[T1]{fontenc}
\usepackage[utf8]{inputenc}
\usepackage{graphicx}
\usepackage{fancyhdr}
\usepackage{float}
\usepackage{xcolor}
\usepackage{authblk}
\usepackage{mathrsfs}
\usepackage{empheq}
\usepackage[hyphens]{url}
\usepackage{hyperref} 
\usepackage[]{breakurl}

\usepackage{geometry}
\DeclareMathOperator{\arctanh}{arctanh}

\begin{document}

\title{Relations for the difference of two dilogarithms}

\author[$\dagger$]{Jean-Christophe {\sc Pain}$^{1,2,}$\footnote{jean-christophe.pain@cea.fr}\\
\small
$^1$CEA, DAM, DIF, F-91297 Arpajon, France\\
$^2$Universit\'e Paris-Saclay, CEA, Laboratoire Mati\`ere en Conditions Extr\^emes,\\ 
91680 Bruy\`eres-le-Ch\^atel, France
}

\maketitle

\begin{abstract}
In this work, we propose a double-series representation of the difference between two dilogarithms with specific arguments. The summation is derived combining a formula we recently found for the so-called Grotendieck-Krivine constant, and an identity obtained by Lima for the difference of two dilogarithms with arguments $\sqrt{2}-1$ and $1-\sqrt{2}$ respectively. We also give an integral representation of the difference considered by Lima, on the basis of a formula published by Brychkov.
\end{abstract}

\section{Introduction}

The dilogarithm (or Spence's function \cite{Watson1937}) \cite{Lewin1958,mathe} is defined as
\begin{equation}
\mathrm{Li}_{2}(z)=\sum_{k=1}^{\infty}\frac{z^{k}}{k^{2}}.
\end{equation}
It can be expressed using the integral representation
\begin{equation}
\mathrm{Li}_2(z)=-\int_0^z\frac{\ln(1-t)}{t}dt.
\end{equation}
The dilogarithm finds many applications in mathematics and theoretical physics \cite{Zagier2007,Fock2009}. The efficient computation of dilogarithms is an important field of investigation \cite{Bailey1997}, as well as the search for special values \cite{Borwein2001} and identities \cite{Kirillov1994,Campbell2021}. There is a number of relations involving at least two dilogarithms. The most commonly used are probably the following ones:
\begin{equation}\label{sum}
\mathrm{Li}_{2}(z)+\mathrm {Li}_{2}(-z)={\frac {1}{2}}\mathrm{Li}_{2}(z^{2}),
\end{equation}
\begin{equation}
\mathrm{Li}_{2}(1-z)+\mathrm {Li}_{2}\left(1-{\frac{1}{z}}\right)=-{\frac{(\ln z)^{2}}{2}},
\end{equation}
\begin{equation}
\mathrm{Li}_{2}(z)+\mathrm {Li}_{2}(1-z)={\frac{\pi^{2}}{6}}-\ln(z)\ln(1-z),
\end{equation}
\begin{equation}
\mathrm{Li}_{2}(-z)-\mathrm {Li}_{2}(1-z)+{\frac {1}{2}}\mathrm {Li}_{2}(1-z^{2})=-{\frac{\pi^{2}}{12}}-\ln (z)\ln(z+1),
\end{equation}
and
\begin{equation}
\mathrm{Li}_{2}(z)+\mathrm {Li}_{2}\left({\frac{1}{z}}\right)=-{\frac{\pi^{2}}{6}}-{\frac{1}{2}}[\ln(-z)]^{2}.
\end{equation}
Ramanujan gave the identities \cite{Berndt1994,Gordon1997}:
\begin{equation}
\mathrm{Li}_2\left(\frac{1}{3}\right)-\frac{1}{6}\mathrm{Li}_2\left(\frac{1}{9}\right)=\frac{\pi^2}{18}-\frac{1}{6}(\ln 3)^2,
\end{equation}
\begin{equation}
\mathrm{Li}_2\left(-\frac{1}{2}\right)+\frac{1}{6}\mathrm{Li}_2\left(\frac{1}{9}\right)=-\frac{\pi^2}{18}+\ln 2\ln 3-\frac{1}{2}(\ln 2)^2-\frac{1}{3}(\ln 3)^2,
\end{equation}
\begin{equation}
\mathrm{Li}_2\left(\frac{1}{4}\right)+\frac{1}{3}\mathrm{Li}_2\left(\frac{1}{9}\right)=\frac{1}{18}\pi^2+2\ln 2\ln 3-2(\ln 2)^2-\frac{2}{3}(\ln 3)^2,	
\end{equation}
\begin{equation}
\mathrm{Li}_2\left(-\frac{1}{3}\right)-\frac{1}{3}\mathrm{Li}_2\left(\frac{1}{9}\right)=-\frac{\pi^2}{18}+\frac{1}{6}(\ln3)^2,
\end{equation}
\begin{equation}
\mathrm{Li}_2\left(-\frac{1}{8}\right)+\mathrm{Li}_2\left(\frac{1}{9}\right)=-\frac{1}{2}\left[\ln\left(\frac{9}{8}\right)\right]^2,	
\end{equation}
and Bailey \emph{et al.} \cite{Bailey1997} showed that
\begin{equation}
\pi^2=36\,\mathrm{Li}_2\left(\frac{1}{2}\right)-36\,\mathrm{Li}_2\left(\frac{1}{4}\right)-12\,\mathrm{Li}_2\left(\frac{1}{8}\right)+6\,\mathrm{Li}_2\left(\frac{1}{64}\right). 
\end{equation}
In the book by Gradshteyn and Ryzhik (see Ref. \cite{Gradshteyn}, formula 6.254, p. 642), one finds
\begin{equation}
\int_0^{\infty}\mathrm{ci}(x)\sin^2x\frac{dx}{x}=\frac{1}{2}\left[\mathrm{Li}_2\left(\frac{1}{2}\right)-\mathrm{Li}_2\left(-\frac{1}{2}\right)\right],    
\end{equation}
where
\begin{equation}
\mathrm{ci}(x)=-\int_x^{\infty}\frac{\cos\xi}{\xi}d\xi.   
\end{equation}
Recently, Lima obtained the following interesting relation \cite{Lima2012}:
\begin{equation}\label{lima}
\mathrm{Li}_2(\sqrt{2}-1)-\mathrm{Li}_2(1-\sqrt{2})=\frac{\pi^2}{8}-\frac{1}{4}\ln^2\left(\sqrt{2}+1\right),
\end{equation}
for which Stewart proposed an alternative derivation in 2022 \cite{Stewart2022}. In section \ref{sec2}, we emphasize the link between Lima's formula and the Grothendieck-Krivine bound, and propose a double nested-sum expansion for the difference $\mathrm{Li}_2(\sqrt{2}-1)-\mathrm{Li}_2(1-\sqrt{2})$. In section \ref{sec3}, we mention an integral representation of the latter specific difference considered by Lima, based on a general expression available in Brychkov's handbook \cite{Brychkov}. We also point out an integral form for the sum $\mathrm{Li}_2(ab)+\mathrm{Li}_2(-ab)$.

\section{A double nested summation for the difference of two dilogarithms}\label{sec2}

We recently found the following identity \cite{Pain2022}
\begin{equation}\label{gro}
\sum_{n=1}^{\infty}\left[\sum_{k=n+1}^{\infty}(-1)^{k}\left(\frac{1}{4k-1}-\frac{1}{4k-3}\right)\right]^2=\frac{\pi}{16}-\frac{1}{4}\ln^2\left(\sqrt{2}+1\right)
\end{equation}
which involves, as Eq. (\ref{lima}), the constant
\begin{equation}
\ln^2\left(\sqrt{2}+1\right)=\frac{\pi^2}{4K_G^2},
\end{equation}
where $K_G$ is sometimes referred to as the Grothendieck-Krivine constant (see Appendix). Combining Eq. (\ref{lima}) and (\ref{gro}) provides the following infinite-sum representation of the difference between the two specific dilogarithms
\begin{empheq}[box=\fbox]{align}
\mathrm{Li}_2(\sqrt{2}-1)-\mathrm{Li}_2(1-\sqrt{2})=\frac{\pi(\pi-1)}{8}+2\sum_{n=1}^{\infty}\left[\sum_{k=n+1}^{\infty}(-1)^{k}\left(\frac{1}{4k-1}-\frac{1}{4k-3}\right)\right]^2.
\end{empheq}

\section{Integral representation}\label{sec3}

Brychkov (Ref. \cite{Brychkov}, formula 2, section 4.1.6, p.155) gives the following results
\begin{equation}\label{rer}
\int_0^{\infty}\sqrt{a^2-x^2}\,\arcsin(bx)dx=\frac{a}{4b}\left\{\frac{(1-a^2b^2)}{2ab}\ln\left(\frac{1+ab}{1-ab}\right)+ab\left[\mathrm{Li}_2(ab)-\mathrm{Li}_2(-ab)\right]-1\right\}    
\end{equation}
with $|\arg(1-a^2b^2)|<\pi$. Equation (\ref{rer}) can be put, in particular using Eq. (\ref{sum}), in the following form:
\begin{eqnarray}
\int_0^{\infty}\sqrt{a^2-x^2}\arcsin(bx)dx&=&\frac{1}{8b^2}\left\{2(1-a^2b^2)\,\arctanh(ab)-2ab\right.\nonumber\\
& &\left.-a^2b^2\left[\mathrm{Li}_2(a^2b^2)-4\,\mathrm{Li}_2(ab)\right]\right\}.   
\end{eqnarray}
In the special case where $a=1$ and $b=\sqrt{2}-1$, one gets
\begin{empheq}[box=\fbox]{align}
\mathrm{Li}_2(\sqrt{2}-1)-\mathrm{Li}_2(1-\sqrt{2})&=4\int_0^1\sqrt{1-x^2}\,\arcsin\left[(\sqrt{2}-1)x\right]dx\nonumber\\
& +(1+\sqrt{2})\left[1+\ln(2-\sqrt{2})-\frac{\ln 2}{2}\right].
\end{empheq} 
It is worth mentioning that Valdebenito obtained the following expressions \cite{Valdebenito2019}:
\begin{equation}
\pi^2=4\left[\ln(1+\sqrt{2})\right]^2-16\int_{\ln(1+\sqrt{2})/2}^{\infty}\ln(\tanh x)dx,
\end{equation}
and
\begin{equation}
\pi^2=-4\left[\ln(1+\sqrt{2})\right]^2-16\int_{0}^{\ln(1+\sqrt{2})/2}\ln(\tanh x)dx,
\end{equation}
from which we find that
\begin{empheq}[box=\fbox]{align}
\mathrm{Li}_2(\sqrt{2}-1)-\mathrm{Li}_2(1-\sqrt{2})=-\frac{1}{2}\int_0^{\ln(1+\sqrt{2})/2}\ln(\tanh x)dx-\frac{3}{2}\int_{\ln(1+\sqrt{2})/2}^{\infty}\ln(\tanh x)dx,
\end{empheq} 
or equivalently
\begin{equation}
\mathrm{Li}_2(\sqrt{2}-1)-\mathrm{Li}_2(1-\sqrt{2})=-\frac{3}{2}\int_0^{\infty}\ln(\tanh x)dx+\int_{0}^{\ln(1+\sqrt{2})/2}\ln(\tanh x)dx.
\end{equation}

We would like to point out, that Brychkov (Ref. \cite{Brychkov}, formula 1, section 4.1.6, p.155) also gives (see also Ref. \cite{Salahuddin2020}):
\begin{equation}
\int_0^{a}\frac{\arcsin[(\sqrt{2}-1)x]}{\sqrt{a^2-x^2}}dx=\frac{1}{2}\left[\mathrm{Li}_2(ab)+\mathrm{Li}_2(-ab)\right],
\end{equation}
with $|\arg(1-a^2b^2)|<\pi$. Using Eq. (\ref{sum}), one can write
\begin{equation}
\int_0^{a}\frac{\arcsin(bx)}{\sqrt{a^2-x^2}}dx=\frac{1}{4}\left[4\,\mathrm{Li}_2(ab)-\mathrm{Li}_2(a^2b^2)\right].
\end{equation}
In the case where $a=1$ and $b=\sqrt{2}-1$, one gets
\begin{equation}
\int_0^{1}\frac{\arcsin[(\sqrt{2}-1)x]}{\sqrt{1-x^2}}dx=\frac{1}{4}\left[4\,\mathrm{Li}_2(\sqrt{2}-1)-\mathrm{Li}_2(3-2\sqrt{2})\right].
\end{equation}
In fact, one has
\begin{equation}
\frac{1}{2}\left[\mathrm{Li}_s(z)-\mathrm{Li}_s(-z)\right]=\chi_s(z),
\end{equation}
where $\chi_s(z)$ is the so-called Legendre chi function, equal to
\begin{equation}
\chi_s(z)=\sum_{k=0}^{\infty}\frac{z^{2k+1}}{(2k+1)^s},
\end{equation}
for $|z|\leq 1$ and $\Re(s)>1$. $\chi_s(z)$ resembles the Dirichlet series for the polylogarithm function $\mathrm{Li}_s(z)$. Nice reviews of the theory of such functions are given by Lewin \cite{Lewin1958,Lewin1981} and Berndt \cite{Berndt1994}. Cvijovi\'c published  integral representations of the Legendre chi functio \cite{Cvijovic2007}, which are thus likely to provide, via $\chi_2(z)$, expressions for $\mathrm{Li}_2(z)-\mathrm{Li}_2(-z)$.

\section{Conclusion}

In this short note, we pointed out an intriguing similarity between a relation obtained by Lima for the difference of two dilogarithms with arguments $\sqrt{2}-1$ and $1-\sqrt{2}$ respectively, and a double-series representation involving the so-called Grotendieck-Krivine constant. The combination of the latter two identities provides an additional relation for the difference of two dilogarithms evaluated for the special values $\sqrt{2}-1$ and $1-\sqrt{2}$. We also provided an integral representation of the difference considered by Lima, using a formula published by Brychkov.

\section*{Appendix: On the Grothendieck-Krivine constant}

Let ${\bf A}$ be an $n\times n$ real square matrix with $n\ge 2$ such that:
\begin{equation}
\left|\sum_{i=1}^n\sum_{j=1}^na_{ij}s_it_j\right|\leq 1 	
\end{equation}
for all real numbers $s_1$, $s_2$, ..., $s_n$ and $t_1$, $t_2$, ..., $t_n$ satisfying $\left|s_i\right|$, $\left|t_j\right|\leq 1$. Then Grothendieck showed that there exists a constant $K_R(n)$ ensuring
\begin{equation}
\left|\sum_{i=1}^n\sum_{j=1}^na_{ij}x_i.y_j\right|\leq K_R(n) 	
\end{equation}
for all vectors $x_1$, $x_2$,..., $x_m$ and $y_1$, $y_2$,..., $y_n$ in a Hilbert space with norms $\left|x_i\right|\leq 1$ and $\left|y_j\right|\leq 1$. The Grothendieck constant is the smallest possible value of $K_R(n)$. As mentioned in the introduction, Krivine postulated that the limit 
\begin{equation}
\lim_{n\rightarrow \infty}K_R(n)	
\end{equation}
is equal to \cite{Krivine1977,Krivine1979}:
\begin{equation}\label{KG2}
K_G=\frac{\pi}{2\ln(1+\sqrt{2})}. 	
\end{equation}

\end{document}